\documentclass[12pt]{article}
\usepackage{amsfonts}
\usepackage{full page,
amssymb, amscd, amsmath,graphicx, %showkeys
}
    \title{{\bf Tensor products of strongly graded vertex algebras
    and their modules}}
    \author{Jinwei Yang}
    \date{}
    
\begin{document}
    \bibliographystyle{alpha}
    \maketitle

    \newtheorem{rema}{Remark}[section]
    \newtheorem{propo}[rema]{Proposition}
    \newtheorem{theo}[rema]{Theorem}
   \newtheorem{defi}[rema]{Definition}
    \newtheorem{lemma}[rema]{Lemma}
    \newtheorem{corol}[rema]{Corollary}
     \newtheorem{exam}[rema]{Example}
\newtheorem{assum}[rema]{Assumption}
     \newtheorem{nota}[rema]{Notation}
        \newcommand{\ba}{\begin{array}}
        \newcommand{\ea}{\end{array}}
        \newcommand{\be}{\begin{equation}}
        \newcommand{\ee}{\end{equation}}
        \newcommand{\bea}{\begin{eqnarray}}
        \newcommand{\eea}{\end{eqnarray}}
        \newcommand{\nno}{\nonumber}
        \newcommand{\nn}{\nonumber\\}
        \newcommand{\lbar}{\bigg\vert}
        \newcommand{\p}{\partial}
        \newcommand{\dps}{\displaystyle}
        \newcommand{\bra}{\langle}
        \newcommand{\ket}{\rangle}
 \newcommand{\res}{\mbox{\rm Res}}
\newcommand{\wt}{\mbox{\rm wt}\;}
\newcommand{\swt}{\mbox{\scriptsize\rm wt}\;}
 \newcommand{\pf}{{\it Proof}\hspace{2ex}}
 \newcommand{\epf}{\hspace{2em}$\square$}
 \newcommand{\epfv}{\hspace{1em}$\square$\vspace{1em}}
        \newcommand{\ob}{{\rm ob}\,}
        \renewcommand{\hom}{{\rm Hom}}
\newcommand{\C}{\mathbb{C}}
\newcommand{\R}{\mathbb{R}}
\newcommand{\Z}{\mathbb{Z}}
\newcommand{\N}{\mathbb{N}}
\newcommand{\A}{\mathcal{A}}
\newcommand{\Y}{\mathcal{Y}}
\newcommand{\comp}{\mathrm{COMP}}
\newcommand{\lgr}{\mathrm{LGR}}

\newcommand{\dlt}[3]{#1 ^{-1}\delta \bigg( \frac{#2 #3 }{#1 }\bigg) }

\newcommand{\dlti}[3]{#1 \delta \bigg( \frac{#2 #3 }{#1 ^{-1}}\bigg) }

 \makeatletter
\newlength{\@pxlwd} \newlength{\@rulewd} \newlength{\@pxlht}
\catcode`.=\active \catcode`B=\active \catcode`:=\active
\catcode`|=\active
\def\sprite#1(#2,#3)[#4,#5]{
   \edef\@sprbox{\expandafter\@cdr\string#1\@nil @box}
   \expandafter\newsavebox\csname\@sprbox\endcsname
   \edef#1{\expandafter\usebox\csname\@sprbox\endcsname}
   \expandafter\setbox\csname\@sprbox\endcsname =\hbox\bgroup
   \vbox\bgroup
  \catcode`.=\active\catcode`B=\active\catcode`:=\active\catcode`|=\active
      \@pxlwd=#4 \divide\@pxlwd by #3 \@rulewd=\@pxlwd
      \@pxlht=#5 \divide\@pxlht by #2
      \def .{\hskip \@pxlwd \ignorespaces}
      \def B{\@ifnextchar B{\advance\@rulewd by \@pxlwd}{\vrule
         height \@pxlht width \@rulewd depth 0 pt \@rulewd=\@pxlwd}}
      \def :{\hbox\bgroup\vrule height \@pxlht width 0pt depth
0pt\ignorespaces}
      \def |{\vrule height \@pxlht width 0pt depth 0pt\egroup
         \prevdepth= -1000 pt}
   }
\def\endsprite{\egroup\egroup}
\catcode`.=12 \catcode`B=11 \catcode`:=12 \catcode`|=12\relax
\makeatother

\def\hboxtr{\FormOfHboxtr} % Only necessary if
%\kern... is wanted
\sprite{\FormOfHboxtr}(25,25)[0.5 em, 1.2 ex] % Resolution ca. 200x340 dpi.

:BBBBBBBBBBBBBBBBBBBBBBBBB | :BB......................B |
:B.B.....................B | :B..B....................B |
:B...B...................B | :B....B..................B |
:B.....B.................B | :B......B................B |
:B.......B...............B | :B........B..............B |
:B.........B.............B | :B..........B............B |
:B...........B...........B | :B............B..........B |
:B.............B.........B | :B..............B........B |
:B...............B.......B | :B................B......B |
:B.................B.....B | :B..................B....B |
:B...................B...B | :B....................B..B |
:B.....................B.B | :B......................BB |
:BBBBBBBBBBBBBBBBBBBBBBBBB |

\endsprite

\def\shboxtr{\FormOfShboxtr} % Only necessary if
%\kern... is wanted
\sprite{\FormOfShboxtr}(25,25)[0.3 em, 0.72 ex] % Resolution ca. 200x340 dp

:BBBBBBBBBBBBBBBBBBBBBBBBB | :BB......................B |
:B.B.....................B | :B..B....................B |
:B...B...................B | :B....B..................B |
:B.....B.................B | :B......B................B |
:B.......B...............B | :B........B..............B |
:B.........B.............B | :B..........B............B |
:B...........B...........B | :B............B..........B |
:B.............B.........B | :B..............B........B |
:B...............B.......B | :B................B......B |
:B.................B.....B | :B..................B....B |
:B...................B...B | :B....................B..B |
:B.....................B.B | :B......................BB |
:BBBBBBBBBBBBBBBBBBBBBBBBB |

\endsprite

\vspace{2em}

%\tableofcontents \vspace{1em}
%\noindent{\large \bf References}\hfill %290

%\newpage

\renewcommand{\theequation}{\thesection.\arabic{equation}}
\renewcommand{\therema}{\thesection.\arabic{rema}}
\setcounter{equation}{0} \setcounter{rema}{0}
\date{}
\maketitle

\begin{abstract}
We study strongly graded vertex algebras and their strongly graded
modules, which are conformal vertex algebras and their modules with a
second, compatible grading by an abelian group satisfying certain
grading restriction conditions. We consider a tensor product of
strongly graded vertex algebras and its tensor product strongly graded
modules. We prove that a tensor product of strongly graded irreducible
modules for a tensor product of strongly graded vertex algebras is
irreducible, and that such irreducible modules, up to equivalence,
exhaust certain naturally defined strongly graded irreducible modules
for a tensor product of strongly graded vertex algebras. We also prove
that certain naturally defined strongly graded modules for the tensor
product strongly graded vertex algebra are completely reducible if and
only if every strongly graded module for each of the tensor product
factors is completely reducible. These results generalize the
corresponding known results for vertex operator algebras and their
modules.
\end{abstract}

\section{Introduction}
We prove that a tensor product of strongly graded irreducible
modules for a tensor product of strongly graded vertex algebras is
irreducible, and that conversely, such irreducible modules, up to
equivalence, exhaust certain naturally defined strongly graded
irreducible modules for a tensor product of strongly graded vertex
algebras. (These terms are defined below.) As a consequence, we
determine all the strongly graded
irreducible modules for the tensor product of the moonshine module
vertex operator algebra $V^{\natural}$ with a vertex algebra
associated with a self-dual even lattice, in particular, the
two-dimensional Lorentzian lattice.

The moonshine conjecture of Conway and Norton in [CN] included the
conjecture that there should exist an infinite-dimensional
representation $V$ of the (not yet constructed) Fischer-Griess
Monster sporadic finite simple group $\mathbb{M}$  such that the
McKay-Thompson series $T_g$ for $g \in \mathbb{M}$ acting on $V$
should have coefficients that are equal to the coefficients of the
$q$-series expansions of certain modular functions. In particular,
this conjecture incorporated the McKay-Thompson conjecture, which asserted that
there should exist a (suitably nontrivial) $\mathbb{Z}$-graded $\mathbb{M}$-module $V =
\coprod_{i \geq -1} V_{-i}$ with graded dimension equal to the
elliptic modular function $j(\tau) - 744 = \sum_{i \geq -1}
c(i)q^i$, where we write $q$ for $e^{2\pi i \tau}$, $\tau$ in the
upper half-plane. Such an $\mathbb{M}$-module, the ``moonshine
module," denoted by $V^{\natural}$, was constructed in \cite{FLM},
and in fact, the construction of \cite{FLM} gave a vertex operator
algebra structure on $V^{\natural}$ equipped with an action of $\mathbb{M}$. In
\cite{FLM}, the authors also gave an explicit formula for the
McKay-Thompson series of any element of the centralizer of an
involution of type $2\mathbf{B}$ of $\mathbb{M}$; the case of the
identity element of $\mathbb{M}$ proved the McKay-Thompson
conjecture.

Borcherds then showed in \cite{B} that the rest of the
McKay-Thompson series for the elements of $\mathbb{M}$ acting on
$V^{\natural}$ are the expected modular functions. He obtained
recursion formulas for the coefficients of McKay-Thompson series for
$V^{\natural}$ from the Euler-Poincar$\acute{{\rm e}}$ identity for
certain homology groups associated with a special Lie algebra, the
``monster Lie
algebra," which he constructed using the tensor product of the
moonshine module vertex operator algebra $V^{\natural}$ and a
natural vertex
algebra associated with the two-dimensional Lorentzian lattice.
The importance of this tensor product vertex algebra motivates the
present paper.

The difference between the terminology ``vertex operator algebra," as
defined in \cite{FLM}, and ``vertex algebra," as defined in
\cite{B}, is that a vertex operator algebra amounts to a vertex
algebra with a conformal vector such that the eigenspaces of the
operator $L(0)$ are all finite dimensional with (integral) eigenvalues
that are truncated from below (cf. \cite{LL}). In \cite{HLZ}, the
authors use a notion of ``conformal vertex algebra," which is a
vertex algebra with a conformal vector and with an $L(0)$-eigenspace
decomposition, and a notion of ``strongly
graded conformal vertex algebra," which is a conformal vertex
algebra with a second, compatible grading by an abelian group
satisfying certain grading restriction conditions.

In a series of papers (\cite{tensorK}--\cite{tensor3},
\cite{tensor4}), the authors developed a tensor product theory for
modules for a vertex operator algebra under suitable conditions. A
structure called ``vertex tensor category structure,'' which is much
richer than braided tensor category structure, has thereby been established
for many important categories of modules for classes of vertex
operator algebras (see \cite{tensorK}). It is expected that a vertex
tensor category together with certain additional structures
determines uniquely (up to isomorphism) a vertex operator algebra
such that the vertex tensor category constructed from a suitable
category of modules for it is equivalent (in the sense of vertex
tensor categories) to the original vertex tensor category. In
\cite{HLZ}, this tensor product theory is generalized to a larger
family of categories of ``strongly graded modules" for a conformal
vertex algebra, under suitably relaxed conditions. We want to
investigate the vertex tensor category in the sense of
\cite{tensorK}, but in the setting of \cite{HLZ}, associated with
the tensor product of the moonshine module vertex operator algebra
$V^{\natural}$ and the vertex algebra associated with the
two-dimensional Lorentzian lattice. The first step in thinking about
this is to determine the irreducible modules for this algebra.

For the vertex operator algebra case, it is proved in \cite{FHL}
that a tensor product module $W_1 \otimes \cdots \otimes W_p$ for a
tensor product vertex algebra $V_1 \otimes \cdots \otimes V_p$
(where $W_i$ is a $V_i$-module) is irreducible if and only if each
$W_i$ is irreducible. The proof uses a version of Schur's Lemma and
also the density theorem \cite{J}. It is also proved in \cite{FHL}
that these irreducible modules $W$ are (up to equivalence) exactly
all the irreducible modules for the tensor product algebra $V_1
\otimes \cdots \otimes V_p$. The proof uses the fact that each
homogeneous subspace of $W$ is finite dimensional. In this paper, we
generalize the arguments in \cite{FHL} to
prove similar, more general results for strongly graded modules for
strongly graded conformal vertex algebras.

For the strongly graded conformal vertex algebra case, the
homogeneous subspaces of a strongly graded module are no longer
finite dimensional. However, by using the fact that each doubly
homogeneous subspace (homogeneous with respect to both gradings) of a
strongly graded conformal vertex algebra is
finite dimensional, we prove a suitable version of Schur's Lemma for
strongly graded modules under the assumption that the abelian group
that gives the second grading of the strongly graded algebra is
countable.

To avoid unwanted flexibility in the second grading such as a
shifting of the grading by an element of the abelian group, we
suppose that the grading abelian groups $A$ for a strongly graded
conformal vertex algebra and $\tilde{A}$ (which includes $A$ as a
subgroup) for its strongly graded modules are always determined by a
vector space, which we typically call $\mathfrak{h}$, consisting of
operators induced by $V$. We call this kind of strongly graded
conformal vertex algebra a ``strongly $(\mathfrak{h}, A)$-graded
conformal vertex algebra" and its strongly graded modules ``strongly
$(\mathfrak{h}, \tilde{A})$-graded modules." Important examples of
strongly $(\mathfrak{h}, A)$-graded conformal vertex algebras and
their strongly $(\mathfrak{h}, \tilde{A})$-graded modules are the
vertex algebras associated with nondegenerate even lattices and
their modules.

For strongly $(\mathfrak{h}_i, \tilde{A_i})$-graded modules $W_i$
for strongly $(\mathfrak{h}_i, A_i)$-graded conformal vertex
algebras $V_i$, we construct a tensor product strongly $(\oplus_{i
= 1}^p \mathfrak{h}_i, \oplus_{i = 1}^p \tilde{A_i})$-graded module
$W_1 \otimes \cdots \otimes W_p$ for the tensor product strongly
graded conformal vertex algebra $V_1 \otimes \cdots \otimes V_p$.
Then we prove that this tensor product module $W_1 \otimes \cdots \otimes
W_p$ is irreducible if and only if each $W_i$ is irreducible, under
the assumption that each grading abelian group $A_i$ for $V_i$ is a
countable group.

To determine all the irreducible strongly graded modules (up to
equivalence) for the tensor product strongly graded conformal vertex
algebra $V_1 \otimes \cdots \otimes V_p$, the main difficulty is
that we need to deal with the second grading by the abelian groups.
For the strongly $(\oplus_{i = 1}^p \mathfrak{h}_i,
\tilde{A})$-graded modules $W$ for the tensor product strongly $(
\oplus_{i = 1}^p \mathfrak{h}_i, \oplus_{i = 1}^p A_i)$-graded
vertex algebra $V_1 \otimes \cdots \otimes V_p$, we assume there is
a decomposition $\tilde{A} = \tilde{A_1}\oplus \cdots \oplus
\tilde{A_p}$, such that $W$ is an $(\mathfrak{h}_i,
\tilde{A_i})$-graded module (that is, a strongly graded module
except for the grading restriction conditions) when viewed as a
$V_i$-module. We call this kind of strongly $(\oplus_{i = 1}^p
\mathfrak{h}_i, \tilde{A})$-graded module a strongly
$((\mathfrak{h}_1, \tilde{A_1})$, \dots, $(\mathfrak{h}_p,
\tilde{A_p}))$-graded module. In the main theorem, we prove that if
such a module is irreducible, then it is a tensor product of
strongly graded irreducible modules. Then, as a corollary of the main theorem, we classify the
strongly graded modules for the tensor
product strongly graded conformal vertex algebra $V^{\natural}
\otimes V_L$, where $L$ is an even lattice, and in particular, where $L$ is the
(self-dual) two-dimensional Lorentzian lattice.

It is proved in \cite{DMZ} that every module for the tensor product
vertex operator algebra $V_1 \otimes \cdots \otimes V_p$ is
completely reducible if and only if every module for each vertex
operator algebra $V_i$ is completely reducible. We also generalize the
argument in \cite{DMZ} to prove a similar result for tensor product
strongly $(\mathfrak{h}, A)$-graded conformal vertex algebras.

This paper is organized as follows: In Section 2, we introduce the
definitions and some basic properties of strongly graded vertex
algebras and their strongly graded modules. Then we construct a tensor
product of strongly graded vertex algebras and its tensor product
strongly graded modules in Section 3. In Section 4, we introduce the
definition of strongly $(\mathfrak{h}, A)$-graded vertex algebra
and strongly $(\mathfrak{h}, \tilde{A})$-graded module. In Section
5, we prove the main theorem, which classifies the irreducible
strongly $((\mathfrak{h}_1, \tilde{A_1})$, \dots,
$(\mathfrak{h}_p, \tilde{A_p}))$-graded $V_1\otimes \cdots \otimes
V_p$-modules. Then we use the main theorem to determine all the
strongly graded modules for $V^{\natural} \otimes V_L$. In Section 6, we consider strongly
graded conformal vertex algebras whose strongly graded modules are
all completely reducible and prove that every strongly
$((\mathfrak{h}_1, \tilde{A_1})$, \dots, $(\mathfrak{h}_p,
\tilde{A_p}))$-graded module for the tensor product strongly graded
algebra $V_1\otimes \cdots \otimes V_p$ is completely reducible if
and only
if every strongly $(\mathfrak{h}_i, \tilde{A_i})$-graded module for
each $V_i$ is completely reducible.\\

\noindent {\bf Acknowledgements} I would like to thank Professor
James Lepowsky for helpful discussions and suggestions. I would also
like to thank Professor Yi-Zhi Huang and Professor Haisheng Li for useful
discussions.

\section{Strongly graded vertex algebras and their modules}

We recall the following four definitions from [HLZ].

\begin{defi}\label{cva}
{\rm A {\it conformal vertex algebra} is a ${\mathbb Z}$-graded
vector space
\begin{equation}\label{Vgrading}
V=\coprod_{n\in{\mathbb Z}} V_{(n)}
\end{equation}
(for $v\in V_{(n)}$, we say the {\it weight} of $v$ is $n$ and we
write $\mbox{wt}\, v=n$) equipped with a linear map $V\otimes V\to
V[[x, x^{-1}]]$, or equivalently,
\begin{eqnarray}\label{YforV}
V&\to&({\rm End}\; V)[[x, x^{-1}]] \nno\\
v&\mapsto& Y(v, x)={\displaystyle \sum_{n\in{\mathbb
Z}}}v_{n}x^{-n-1} \;\;( \mbox{where }v_{n}\in{\rm End} \;V),
\end{eqnarray}
$Y(v, x)$ denoting the {\it vertex operator associated with} $v$,
and equipped also with two distinguished vectors ${\bf 1}\in
V_{(0)}$ (the {\it vacuum vector}) and $\omega\in V_{(2)}$ (the {\it
conformal vector}), satisfying the following conditions for $u,v \in
V$: the {\it lower truncation condition}:
\begin{equation}\label{ltc}
u_{n}v=0\;\;\mbox{ for }n\mbox{ sufficiently large}
\end{equation}
(or equivalently, $Y(u, x)v\in V((x))$); the {\it vacuum property}:
\begin{equation}\label{1left}
Y({\bf 1}, x)=1_V;
\end{equation}
the {\it creation property}:
\begin{equation}\label{1right}
Y(v, x){\bf 1} \in V[[x]]\;\;\mbox{ and }\;\lim_{x\rightarrow 0}Y(v,
x){\bf 1}=v
\end{equation}
(that is, $Y(v, x){\bf 1}$ involves only nonnegative integral powers
of $x$ and the constant term is $v$); the {\it Jacobi identity} (the
main axiom):
\begin{eqnarray}
&x_0^{-1}\delta \bigg({\displaystyle\frac{x_1-x_2}{x_0}}\bigg)Y(u,
x_1)Y(v, x_2)-x_0^{-1} \delta
\bigg({\displaystyle\frac{x_2-x_1}{-x_0}}\bigg)Y(v, x_2)Y(u,
x_1)&\nno \\ &=x_2^{-1} \delta
\bigg({\displaystyle\frac{x_1-x_0}{x_2}}\bigg)Y(Y(u, x_0)v,
x_2)&\label{Jacobi}
\end{eqnarray}
(note that when each expression in (\ref{Jacobi}) is applied to any
element of $V$, the coefficient of each monomial in the formal
variables is a finite sum; on the right-hand side, the notation
$Y(\cdot, x_2)$ is understood to be extended in the obvious way to
$V[[x_0, x^{-1}_0]]$); the {\em Virasoro algebra relations}:
\begin{equation}\label{vir1}
[L(m), L(n)]=(m-n)L(m+n)+{\displaystyle\frac{1}{12}}
(m^3-m)\delta_{n+m,0}c
\end{equation}
for $m, n \in {\mathbb Z}$, where
\begin{equation}\label{vir2}
L(n)=\omega _{n+1}\;\; \mbox{ for } \;  n\in{\mathbb Z}, \;\;
\mbox{i.e.},\;\;Y(\omega, x)=\sum_{n\in{\mathbb Z}}L(n)x^{-n-2},
\end{equation}
\begin{equation}\label{vir3}
c\in {\mathbb C}
\end{equation}
(the {\it central charge} or {\it rank} of $V$);
\begin{equation}\label{L-1derivativeproperty}
{\displaystyle \frac{d}{dx}}Y(v, x)=Y(L(-1)v, x)
\end{equation}
(the {\it  $L(-1)$-derivative property}); and
\begin{equation}\label{L0gradingproperty}
L(0)v=nv=(\wt v)v \;\; \mbox{ for }\; n\in {\mathbb Z}\; \mbox{ and
}\; v\in V_{(n)}.
\end{equation}
}
\end{defi}

This completes the definition of the notion of conformal vertex
algebra.  We will denote such a conformal vertex algebra by
$(V,Y,{\bf 1},\omega)$. \\

\begin{defi}\label{cvamodule}
{\rm  Given a conformal vertex algebra $(V,Y,{\bf 1},\omega)$,  a
{\it module} for $V$ is a ${\mathbb C}$-graded vector space
\begin{equation}\label{Wgrading}
W=\coprod_{n\in{\mathbb C}} W_{(n)}
\end{equation}
(graded by {\it weights}) equipped with a linear map $V\otimes W
\rightarrow W[[x,x^{-1}]]$, or equivalently,
\begin{eqnarray}\label{YforW}
V &\rightarrow & (\mbox{End}\ W)[[x,x^{-1}]] \nno\\
v&\mapsto & Y(v,x) =\sum_{n\in {\mathbb Z}}v_nx^{-n-1}\;\;\;
(\mbox{where}\;\; v_n \in \mbox{End}\ W)
\end{eqnarray}
(note that the sum is over ${\mathbb Z}$, not ${\mathbb C}$),
$Y(v,x)$ denoting the {\it vertex operator on $W$ associated with
$v$}, such that all the defining properties of a conformal vertex
algebra that make sense hold.  That is, the following conditions are
satisfied: the lower truncation condition: for $v \in V$ and $w \in
W$,
\begin{equation}\label{ltc-w}
v_nw = 0 \;\;\mbox{ for }\;n\;\mbox{ sufficiently large}
\end{equation}
(or equivalently, $Y(v, x)w\in W((x))$); the vacuum property:
\begin{equation}\label{m-1left}
Y(\mbox{\bf 1},x) = 1_W;
\end{equation}
the Jacobi identity for vertex operators on $W$: for $u, v \in V$,
\begin{eqnarray}\label{m-Jacobi}
&{\dps x^{-1}_0\delta \bigg( {x_1-x_2\over x_0}\bigg)
Y(u,x_1)Y(v,x_2) - x^{-1}_0\delta \bigg( {x_2-x_1\over -x_0}\bigg)
Y(v,x_2)Y(u,x_1)}&\nno\\
&{\dps = x^{-1}_2\delta \bigg( {x_1-x_0\over x_2}\bigg)
Y(Y(u,x_0)v,x_2)}
\end{eqnarray}
(note that on the right-hand side, $Y(u,x_0)$ is the operator on $V$
associated with $u$); the Virasoro algebra relations on $W$ with
scalar $c$ equal to the central charge of $V$:
\begin{equation}\label{m-vir1}
[L(m), L(n)]=(m-n)L(m+n)+{\displaystyle\frac{1}{12}}
(m^3-m)\delta_{n+m,0}c
\end{equation}
for $m,n \in {\mathbb Z}$, where
\begin{equation}\label{m-vir2}
L(n)=\omega _{n+1}\;\; \mbox{ for }n\in{\mathbb Z}, \;\;{\rm
i.e.},\;\;Y(\omega, x)=\sum_{n\in{\mathbb Z}}L(n)x^{-n-2};
\end{equation}
\begin{equation}\label{L-1}
\displaystyle \frac{d}{dx}Y(v, x)=Y(L(-1)v, x)
\end{equation}
(the $L(-1)$-derivative property); and
\begin{equation}\label{wl0}
(L(0)-n)w=0\;\;\mbox{ for }\;n\in {\mathbb C}\;\mbox{ and }\;w\in
W_{(n)},
\end{equation}
where $n = {\rm wt}\ w$.}
\end{defi}

This completes the definition of the notion of module for a
conformal vertex algebra.

\begin{defi}\label{def:dgv}
{\rm Let $A$ be an abelian group.  A conformal vertex algebra
\[
V=\coprod_{n\in {\mathbb Z}} V_{(n)}
\]
is said to be {\em strongly graded with respect to $A$} (or {\em
strongly $A$-graded}, or just {\em strongly graded} if the abelian
group $A$ is understood) if it is equipped with a second gradation, by $A$,
\[
V=\coprod _{\alpha \in A} V^{(\alpha)},
\]
such that the following conditions are satisfied: the two gradations
are compatible, that is,
\[
V^{(\alpha)}=\coprod_{n\in {\mathbb Z}} V^{(\alpha)}_{(n)} \;\;
(\mbox{where}\;V^{(\alpha)}_{(n)}=V_{(n)}\cap V^{(\alpha)})\;
\mbox{ for any }\;\alpha \in A;
\]
for any $\alpha,\beta\in A$ and $n\in {\mathbb Z}$,
\begin{eqnarray}
&V^{(\alpha)}_{(n)}=0\;\;\mbox{ for }\;n\;\mbox{ sufficiently
negative};&\label{dua:ltc}\\
&\dim V^{(\alpha)}_{(n)} <\infty;&\label{dua:fin}\\
&{\bf 1}\in V^{(0)}_{(0)};&\\
&\omega\in V^{(0)}_{(2)};&\\
&v_l V^{(\beta)} \subset V^{(\alpha+\beta)}\;\; \mbox{ for any
}\;v\in V^{(\alpha)},\;l\in {\mathbb Z}.&\label{v_l-A}
\end{eqnarray} }
\end{defi}

This completes the definition of the notion of strongly $A$-graded
conformal vertex algebra.\\

For modules for a strongly graded algebra we will also have a second
grading by an abelian group, and it is natural to allow this group
to be larger than the second grading group $A$ for the algebra.
(Note that this already occurs for the {\em first} grading group,
which is ${\mathbb Z}$ for algebras and ${\mathbb C}$ for modules.)

\begin{defi}\label{def:dgw}{\rm
Let $A$ be an abelian group and $V$ a strongly $A$-graded conformal
vertex algebra. Let $\tilde A$ be an abelian group containing $A$ as
a subgroup. A $V$-module
\[
W=\coprod_{n\in{\mathbb C}} W_{(n)}
\]
is said to be {\em strongly graded with respect to $\tilde A$} (or
{\em strongly $\tilde A$-graded}, or just {\em strongly graded} if
the abelian group $\tilde A$ is understood) if it is equipped with a
second gradation, by $\tilde A$,
\begin{equation}\label{2ndgrd}
W=\coprod _{\beta \in \tilde A} W^{(\beta)},
\end{equation}
such that the following conditions are satisfied: the two gradations
are compatible, that is, for any $\beta \in \tilde A$,
\[
W^{(\beta)}=\coprod_{n\in {\mathbb C}} W^{(\beta)}_{(n)}
\;\;(\mbox{where }\; W^{(\beta)}_{(n)}=W_{(n)}\cap W^{(\beta)})
\]
for any $\alpha\in A$, $\beta\in \tilde A$ and $n\in {\mathbb C}$,
\begin{eqnarray}
&W^{(\beta)}_{(n+k)}=0
\mbox{ for }\;k\in {\mathbb Z}\;\mbox{
sufficiently
negative};&\label{set:dmltc}\\
&\dim W^{(\beta)}_{(n)} <\infty
&\label{set:dmfin}\\
&v_l W^{(\beta)} \subset W^{(\alpha+\beta)}\;\;\mbox{ for any
}\;v\in V^{(\alpha)},\;l\in {\mathbb Z}.&\label{m-v_l-A}
\end{eqnarray} }
\end{defi}

This completes the definition of the notion of strongly
$\tilde{A}$-graded module for a strongly $A$-graded conformal vertex
algebra.

\begin{rema}\label{may shift}{\rm It is always possible that there are
different gradings on $W$
by $\tilde{A}$, such as by shifting by an element in $\tilde{A}$.
However, in this paper, we shall fix one particular
$\tilde{A}$-grading on $W$.}
\end{rema}

In order to study strongly graded $V$-modules for tensor product
algebras, we shall need the following generalization:

\begin{defi}\label{doubly graded}{\rm In the setting of Definition
\ref{def:dgw} (the definition of ``strongly graded module"), a $V$-module
(not necessarily strongly graded, of course) is {\it doubly graded
with respect to $\tilde{A}$} if it satisfies all the conditions in
Definition \ref{def:dgw} except perhaps for (\ref{set:dmltc}) and
(\ref{set:dmfin}).}\end{defi}

\begin{exam}\label{vertex operator algebra}{\rm Note that the notion
of conformal vertex algebra strongly
graded with respect to the trivial group is exactly the notion of
vertex operator algebra. Let $V$ be a vertex operator algebra,
viewed (equivalently) as a conformal vertex algebra strongly graded
with respect to the trivial group.  Then the $V$-modules that are
strongly graded with respect to the trivial group (in the sense of
Definition \ref{def:dgw}) are exactly the ($\C$-graded) modules for
$V$ as a vertex operator algebra, with the grading restrictions as
follows: For $n \in \C$,
\begin{equation}\label{Wn+k=0}
W_{(n+k)}=0 \;\; \mbox { for }\;k\in {\mathbb Z}\;\mbox{
sufficiently negative}
\end{equation}
and
\begin{equation}\label{dimWnfinite}
\dim W_{(n)} <\infty.
\end{equation}}
\end{exam}

\begin{exam}\label{lattice vertex algebra}{\rm
An important source of examples of strongly graded conformal vertex
algebras and modules comes {}from the vertex algebras and modules
associated with even lattices. We recall the following construction
from \cite{FLM}. Let $L$ be an even lattice, i.e., a
finite-rank free abelian group equipped with a nondegenerate
symmetric bilinear form $\langle\cdot,\cdot\rangle$, not necessarily
positive definite, such that $\langle \alpha, \alpha\rangle\in 2{\mathbb Z}$
for all $\alpha\in L$. Let $\mathfrak{h} = L\otimes_{\mathbb{Z}}
\mathbb{C}$. Then $\mathfrak{h}$ is a vector space with a nonsingular
bilinear form $\langle \cdot, \cdot \rangle$, extended from $L$. We form a Heisenberg
algebra
\[\widehat{\mathfrak{h}}_{\mathbb{Z}} = \coprod_{n \in \mathbb{Z},\ n \neq 0}
 \mathfrak{h} \otimes t^n \oplus \mathbb{C}c.\] Let $(\widehat{L},\bar{}\ )$
 be a central extension of $L$ by a finite cyclic group $\langle \kappa\ |\ \kappa^s =
 1\rangle$. Fix a primitive $s$th root of unity, say $\omega$, and define
 the faithful character \[\chi : \langle \kappa \rangle \rightarrow
 \mathbb{C}^{*}\] by the condition \[\chi(\kappa) = \omega.\] Denote by
 $\mathbb{C}_{\chi}$ the one-dimensional space $\mathbb{C}$ viewed as a $\langle \kappa
 \rangle$-module on which $\langle \kappa \rangle$ acts according to
 $\chi$: \[\kappa \cdot 1 = \omega,\] and denote by $\mathbb{C}\{L\}$ the induced
 $\widehat{L}$-module \[\mathbb{C}\{L\} = {\rm
 Ind}^{\widehat{L}}_{\langle \kappa \rangle} \mathbb{C}_{\chi} =
 \mathbb{C}[\widehat{L}]\otimes_{\mathbb{C}[\langle \kappa \rangle]} \mathbb{C}_{\chi}
 .\] Then \[V_L = S(\widehat{\mathfrak{h}}_{\mathbb{Z}}^{-}) \otimes \mathbb{C}\{L\}\] has a
 natural structure of conformal vertex algebra; see \cite{B} and Chapter 8 of
\cite{FLM}. For $\alpha \in L$, choose an $a \in \widehat{L}$ such that
$\bar{a} = \alpha$. Define
\[\iota(a) = a \otimes 1 \in \mathbb{C}\{L\}\] and
\[V_L^{(\alpha)} = {\rm span}\ \{h_1(-n_1)\cdots h_k(-n_k)\otimes
\iota(a)\},\]where $h_1, \dots, h_k \in \mathfrak{h}$,
$n_1, \dots, n_k > 0$, and where $h(n)$ is the operator associated
with $h \otimes t^n$ via the $\hat{\mathfrak{h}}_{\Z}$-module
structure of $V_L$. Then $V_L$ is equipped with a natural second
grading
given by $L$ itself. Also for $n \in \mathbb{Z}$, we have
\[(V_L)_{(n)}^{(\alpha)} = {\rm span}\ \{h_1(-n_1)\cdots
h_k(-n_k)\otimes \iota(a)|\ \bar{a} = \alpha, \sum_{i = 1}^k n_i +
\frac{1}{2}\langle \alpha, \alpha\rangle = n\},\] making $V_L$ a
strongly $L$-graded conformal vertex algebra in the sense of
Definition \ref{def:dgv}. When the form
$\langle\cdot,\cdot\rangle$ on $L$ is also positive definite, then
$V_L$ is a vertex operator algebra, that is, as in Example \ref{vertex
operator algebra}, $V_L$ is a strongly
$A$-graded conformal vertex algebra for $A$ the trivial group. In
general, a conformal vertex algebra may be strongly graded for several
choices of $A$.

Any sublattice $M$ of the ``dual lattice''
$L^{\circ}$ of $L$ containing $L$ gives rise to a strongly
$M$-graded module for the strongly $L$-graded conformal vertex
algebra (see Chapter 8 of \cite{FLM}; cf. \cite{LL}). In fact, any
irreducible $V_L$-module is equivalent to a $V_L$-module of the form
$V_{L + \beta} \subset V_{L^{\circ}}$ for some
$\beta \in L^{\circ}$ and any $V_L$-module $W$ is equivalent to a
direct sum of irreducible $V_L$-modules, i.e., \[W = \coprod_{\gamma_i
\in L^{\circ},\ i = 1, \dots, n} V_{\gamma_i + L},\] where
$\gamma_i$'s are arbitrary elements of $L^{\circ}$, and $n \in
\mathbb{N}$ (see \cite{D1}, \cite{DLM}; cf. \cite{LL}). In general, a
module for a strongly graded vertex algebra may be strongly graded for
several choices of $\tilde{A}$.}
\end{exam}

\begin{nota}\label{VA}{\rm In the remainder of this section, without further assumption,
we will let $A$ be an abelian group and $V$ be a strongly $A$-graded conformal vertex
algebra. Also, we will let $\tilde{A}$ be an abelian group
containing $A$ and $W$ be a doubly graded $V$-module with respect to
$\tilde{A}$. When we need $W$ to be strongly graded, we will say it
explicitly.}\end{nota}

\begin{defi}\label{DH}{\rm
The subspaces $V_{(n)}^{(\alpha)}$ for $n \in \mathbb{Z}$, $\alpha
\in A$ in Definition \ref{doubly graded} are called the {\it doubly
homogeneous subspaces} of $V$. The elements in $V_{(n)}^{(\alpha)}$
are called {\it doubly homogeneous} elements. Similar definitions
can be used for $W^{(\beta)}_{(n)}$ in the module $W$.}
\end{defi}

\begin{nota}\label{A-wt}{\rm Let $v$ be a doubly homogeneous element
of $V$. Let wt $v_n$, $n \in \mathbb{Z}$, refer to the weight of
$v_n$ as an operator acting on $W$, and let $A$-wt $v_n$ refer to
the $A$-weight of $v_n$ on $W$.}\end{nota}

\begin{lemma}\label{b}Let $v \in V^{(\alpha)}_{(n)}$, for $n \in \mathbb{Z}$,
$\alpha \in A$. Then for $m \in \mathbb{Z}$, {\rm wt}$\ v_m$ = $n -
m -1$ and $A$-{\rm wt}$\ v_m = \alpha$.\end{lemma}
\pf The first equation is standard from the theory of graded conformal
vertex algebras and the second follows easily from the
definitions. \epfv

\begin{defi}\label{A(V)}{\rm The algebra $A(V; W)$ {\it associated
with $V$ and $W$} is defined to be
the algebra of operators on $W$ induced by $V$, i.e., the algebra
generated by the set \[ \{ v_n\ |\ v \in V, \ n \in
\mathbb{Z}\}.\]
For a subspace $V^{'}$ of $V$, we use $A(V^{'}; W)$ to denote the
subalgebra of $A(V; W)$ generated by the set
\[ \{ v_n\ |\ v \in V^{'}, \ n \in
\mathbb{Z}\}.\] For a subspace $W^{'}$ of $W$, we use $A(V; W^{'})$
to denote the subalgebra of $A(V; W)$ preserving $W^{'}$. Similarly for $V^{'}$
and $W^{'}$, we use $A(V^{'}; W^{'})$ to denote the subalgebra of
$A(V; W)$ generated by the operators on $W^{'}$ induced by
$V^{'}$.}\end{defi}

\begin{rema}{\rm When $W^{'}$ is a submodule of $W$, there are two
possible definitions for $A(V; W^{'})$ in Definition \ref{A(V)}. One
is as an algebra associated with $V$ and $W^{'}$, the other is as a
subalgebra of $A(V; W)$. But it does not matter because they are
both algebras of operators on $W^{'}$ generated by the set
\[ \{ v_n\ |\ v \in V, \ n \in \mathbb{Z}\}.\]
Similar comments hold for $V^{'}$ a subalgebra of $V$.}\end{rema}

The following lemma follows easily from Lemma \ref{b}:
\begin{lemma}\label{double graded}The algebra $A(V; W)$ is
doubly graded by $\mathbb{Z}$ and $A$.
Moreover for $n \in \mathbb{Z}$, \begin{equation*}\begin{split}A(V;
W)_{(n)} = {\rm span}\ \{(v_1)_{j_1}\cdots (v_m)_{j_m}\ |\ \sum_{i =
1}^m\ {\rm wt}\ (v_i)_{j_i} = n, \\{\rm where}\ m \in \mathbb{N},\
v_i \in V,\ j_i \in \mathbb{Z},\ {\rm for}\ i = 1, \dots,
m\}\end{split}\end{equation*} and for $\alpha \in A$,
\begin{equation*}\begin{split}A(V; W)^{(\alpha)} = {\rm span}\ \{(v_1)_{j_1}\cdots
(v_m)_{j_m}\ |\ \sum_{i = 1}^m A\mbox{-}{\rm wt}\ (v_i)_{j_i} =
\alpha,\\{\rm where}\ m \in \mathbb{N},\ v_i \in V,\ j_i \in
\mathbb{Z},\ {\rm for}\ i = 1, \dots, m\}.\end{split}\end{equation*}\end{lemma}

\begin{propo}\label{c}
Let $W$ be an irreducible doubly graded $V$-module with respect to
$\tilde{A}$. Then we have the following results:

\begin{enumerate}
\item[(a)]Each weight subspace $W_{(h)}$ ($h \in \mathbb{C}$) is
irreducible under the algebra $A(V; W_{(h)})$.

\item[(b)]Each $\tilde{A}$-homogeneous subspace $W^{(\beta)}\ (\beta \in
\tilde{A})$ is irreducible under the algebra $A(V; W^{(\beta)})$.

\item[(c)]Each doubly homogeneous subspace $W_{(h)}^{(\beta)}$ ($h \in
\mathbb{C}$, $\beta \in \tilde{A}$) is irreducible under the algebra
$A(V; W_{(h)}^{(\beta)})$.
\end{enumerate}
\end{propo}

$Proof$. We only prove statement (a), the proofs of statements (b)
and (c) being similar. If $W_{(h)}$ is not irreducible, we can find
a nontrivial proper submodule $U$ of $W_{(h)}$ under the algebra
$A(V; W_{(h)})$. This submodule cannot generate all $W$
under the action by the algebra $A(V; W)$, since by Lemma \ref{double graded}, \[A(V; W)U =
\coprod_{n \in \mathbb{Z}} A(V; W)_{(n)}U \subset U \oplus
\coprod_{m \in \mathbb{Z}, m \neq h} W_{(m)}.
\]This contradicts the irreducibility of $W$.\epfv

\begin{rema}\label{congruent}{\rm
A $V$-module $W$ decomposes into
submodules corresponding to the congruence classes of its weights
modulo $\mathbb{Z}$: For $\mu \in \mathbb{C}/\mathbb{Z}$, let
\begin{equation}
W_{(\mu)} = \coprod_{\bar n = \mu} W_{(n)},
\end{equation}
where $\bar n$ denotes the equivalence class of $n \in \mathbb{C}$
in $\mathbb{C}/\mathbb{Z}$.  Then
\begin{equation}
W = \coprod_{\mu \in \mathbb{C}/\mathbb{Z}} W_{(\mu)}
\end{equation}
and each $W_{(\mu)}$ is a $V$-submodule of $W$.  Thus if a module
$W$ is indecomposable (in particular, if it is irreducible), then
all complex numbers $n$ for which $W_{(n)}\neq 0$ are congruent
modulo ${\mathbb Z}$ to each other. }
\end{rema}

\begin{defi}\label{5}{\rm Let $W_1$ and $W_2$ be doubly
graded $V$-modules with respect to $\tilde{A}$. A
{\it module homomorphism} from $W_1$ to $W_2$ is a linear map $\psi$
such that
\[ \psi(Y(v,x)w)=Y(v,x)\psi(w)\ {\rm for}\ v\in V,\ w\in W_1,\] and
such that $\psi$ preserves the grading by $\tilde{A}$. An {\it
isomorphism} is a bijective homomorphism. An {\it endomorphism} is a
homomorphism from $W$ to itself, we denote the endomorphism ring by
End$_{V}^{\tilde{A}}(W)$.}\end{defi}

\begin{rema}\label{a}{\rm Suppose $V$, $W_1$, $W_2$, $\psi$ are as in Definition
\ref{5}. Then $\psi$ is compatible with both gradings:
\[\psi((W_1)_{(h)}^{(\beta)}) \subset (W_2)_{(h)}^{(\beta)},\ h \in \mathbb{C},
\]
because $\psi$ commutes with $L(0)$ (see Section 4.5 of \cite{LL}),
and because $\psi$ preserves the grading by $\tilde{A}$.}\end{rema}

\begin{rema}\label{commuting ring}{\rm The endomorphism ring
End$_{V}^{\tilde{A}}(W)$ is a subring of the commuting ring
\[ {\rm End}_V(W): = \{{\rm linear\ maps}\ \psi: W \rightarrow W |\
\psi(Y(v,x)w) = Y(v,x)\psi(w), {\rm for}\ v \in V, w \in W
\}.\]}\end{rema}

\begin{propo}\label{xx}Suppose $W$ is an irreducible strongly
$\tilde{A}$-graded $V$-module. Then
${\rm End}^{\tilde{A}}_V (W)$ = $\mathbb{C}$.\end{propo}

$Proof$. For any $\lambda \in \mathbb{C}$, $\psi \in$
End$^{\tilde{A}}_V (W)$, let $W_{\lambda}^{\psi}$ be the
$\lambda$-eigenspace of $\psi$. Then $W_{\lambda}^{\psi}$ is a
$V$-submodule of $W$. Because $W$ is irreducible,
$W_{\lambda}^{\psi}$ = 0 or $W$. We still need to show
$W_{\lambda}^{\psi} \neq 0$, for some $\lambda \in \mathbb{C}$.\\

Choose $h\in \mathbb{C}$, $\beta \in \tilde{A}$ such that
$W_{(h)}^{(\beta)} \neq 0$. Then by Remark \ref{a}, $\psi$ preserves
$W_{(h)}^{(\beta)}$. Since ${\rm dim}\ W_{(h)}^{(\beta)} < \infty$
and we are working over $\mathbb{C}$, $\psi$ has an eigenvector in
$W_{(h)}^{(\beta)}$. Therefore $W_{\lambda}^{\psi}
\neq$ 0 for some $\lambda \in \mathbb{C}$.\epf \\

\begin{propo}\label{x}Suppose $A$ is a countable abelian group. Then {\rm End}$_V(W)$ =
$\mathbb{C}$.\end{propo}

$Proof$. From Definition \ref{def:dgv}, $V_{(n)} = \coprod_{\alpha
\in A} V_{(n)}^{(\alpha)}$, where each doubly homogeneous subspace
$V^{(\alpha)}_{(n)}$ has finite dimension. Since $A$ is a countable
group, there are countably many such doubly homogeneous subspaces
$V^{(\alpha)}_{(n)}$, and hence $V$ has countable dimension. Since
$W$ is irreducible, from Proposition 4.5.6 of \cite{LL}, we know
\[W = {\rm span} \{ v_n w\ |\ v \in V, n \in \mathbb{Z}\},\] for any
nonzero element $w$ in $W$. Since $V$ has countable dimension, so
does $W$. Then the result follows from Dixmier's Lemma, which says
that if $S$ is an irreducible set of operators on a vector space $W$
of countable dimension over $\mathbb{C}$, then the commuting ring of
$S$ on $W$ consists of the scalars (cf.\ Lemma 2.2 in \cite{L}, and
\cite{W}, p.11), where we take $S$ to be $A(V; W)$.
\epf \\

\section{Tensor product of strongly graded vertex algebras and their modules}
In this section, we are going to introduce the notion of tensor
product of finitely many strongly graded conformal vertex algebras
and their modules.

Let $A_1, \dots, A_p$ be abelian groups, and let $V_1, \dots, V_p$
be strongly $A_1, \dots, A_p$-graded conformal vertex algebras with
conformal vectors
$\omega^1, \dots, \omega^p$, respectively.\\

Let
\[A = A_1 \oplus \cdots \oplus A_p.\] Then the vector space
\[V = V_1\otimes \cdots \otimes V_p \] becomes a strongly $A$-graded conformal vertex algebra,
which we shall call the {\it tensor product strongly A-graded
conformal vertex algebra}, with the following structure:
\[
Y(v^{(1)}\otimes \cdots \otimes v^{(p)}, x) = Y(v^{(1)},x)\otimes
\cdots \otimes Y(v^{(p)}, x)\] for $v^{(i)} \in V_i$ and the vacuum
vector is
\[
\bf 1 = 1\otimes \cdots \otimes 1.
\]
(Here we use the notation ${\bf 1}$ for the vacuum vectors of $V$
and each $V_i$.) The conformal vector is
\[ \omega = \omega^1 \otimes {\bf 1 \otimes \cdots \otimes 1 +
\cdots + 1 \otimes \cdots \otimes 1 \otimes} \omega^p.
\]
Then
\[L(n) = L_1(n)\otimes 1 \otimes \cdots \otimes 1 + \cdots + 1\otimes
\cdots \otimes 1 \otimes L_p(n)
\]
or $n \in \mathbb{Z}$. (Here we use the notation $L_i(n)$ for the
operators on $V_i$
associated with $\omega^i$, $i = 1, \dots,
p$.) The $A$-grading of $V$ is given by \[V = \coprod_{\alpha \in A}
V^{(\alpha)},\] with \[ V^{(\alpha)} =
 V_1^{(\alpha_1)}\otimes
\cdots \otimes V_p^{(\alpha_p)},
\]
where $\alpha_i \in A_i$, $i = 1, \dots, p$, are such that $\alpha_1+
\cdots + \alpha_p = \alpha$. The $\mathbb{Z}$-grading of $V$ is
given by
\[V = \coprod_{n \in
\mathbb{Z}} V_{(n)},\] where \[V_{(n)} = \coprod_{n_1 + \cdots + n_p
= n}(V_1)_{(n_1)}\otimes \cdots \otimes (V_p)_{(n_p)}.\] (It follows
that the $\mathbb{Z}$-grading is given by $L(0)$ defined above.)

\begin{propo}\label{V} The tensor product of finitely many strongly
graded conformal vertex algebras is
a strongly graded conformal vertex algebra whose central charge is
the sum of the central charges of the tensor factors.
\end{propo}

{\it Proof.} The grading restrictions (\ref{dua:ltc}) and
(\ref{dua:fin}) clearly hold. The Jacobi identity follows from the
weak commutativity and weak associativity properties, as in Section
3.4 of [LL]. \epf

\begin{nota}{\rm For each $i = 1, \dots,
p$, we identify $V_i$ with the subspace $1 \otimes \cdots \otimes 1
\otimes V_i \otimes 1 \otimes \cdots \otimes 1$ of $V$. The strongly
graded conformal vertex algebra $V_i$ is a vertex subalgebra of $V$.
However, it is not a conformal vertex subalgebra of $V$ because the
conformal vector of $V$ and $V_i$ do not match.}
\end{nota}

\begin{rema}\label{comm}{\rm From the definition of tensor product strongly graded conformal
vertex algebra, we see that \[Y((1 \otimes \cdots \otimes 1 \otimes
v^{(i)} \otimes 1 \otimes \cdots \otimes 1, x) = 1_{V_1} \otimes
\cdots \otimes 1_{V_{i-1}} \otimes Y(v^{(i)}, x) \otimes 1_{V_{i+1}}
\otimes \cdots \otimes 1_{V_p},\] for $v^{(i)} \in V_i$. In
particular, we have
\[ [Y(V_i, x_1), Y(V_j, x_2)] = 0,
\]for $i, j = 1, \dots, p$ and $i\neq j$.}\end{rema}

\begin{lemma}\label{e}For all $n \in$ $\mathbb{Z}$, $(v^{(1)} \otimes
\cdots \otimes v^{(p)})_n$ can be expressed as a linear combination,
finite on any given vector, of operators of the form $(v^{(1)}
\otimes 1 \otimes \cdots \otimes 1)_{i_1} \cdots (1 \otimes \cdots
\otimes 1 \otimes v^{(p)})_{i_p}$.
\end{lemma}

{\it Proof.} We prove the result as in [FHL] by induction.
When $p = 2$, taking Res$_{x_1}$ and the constant term in $x_0$ of
the Jacobi identity, we find that
\begin{equation*}\begin{split}
Y(v^{(1)}\otimes v^{(2)}, x_2) = &\ {\rm Res}_{x_0} x_0^{-1}
Y(Y(v^{(1)} \otimes 1, x_0)(1 \otimes v^{(2)}), x_2)\\
= &\ {\rm Res}_{x_1} (x_1 - x_2)^{-1} Y(v^{(1)}\otimes 1,
x_1)Y(1 \otimes v^{(2)}, x_2)\\
&- \ {\rm Res}_{x_1} (- x_2 + x_1)^{-1} Y(1 \otimes v^{(2)},
x_2)Y(v^{(1)}\otimes 1, x_1), \end{split}\end{equation*} so that for
all $n \in \mathbb{Z}$, $(v^{(1)} \otimes v^{(2)})_n$
can be expressed as a linear
 combination, finite on
 any given vector, of operators of the form $(v^{(1)}
\otimes 1)_{n_1} (1 \otimes v^{(2)})_{n_2}$.(Note that we don't need
operators of the form $(1 \otimes v^{(2)})_{n_2}(v^{(1)} \otimes
1)_{n_1}$ because of the Remark \ref{comm}.)

For general $p$, taking Res$_{x_1}$ and the constant term in $x_0$
of the Jacobi identity, we have
\begin{equation*}\begin{split}
&Y(v^{(1)} \otimes \cdots \otimes v^{(p)}, x_2)\\
 = &\ {\rm
Res}_{x_0} x_0^{-1} Y(Y(v^{(1)}\otimes \cdots \otimes v^{(p-1)}
\otimes 1,
x_0)(1 \otimes \cdots \otimes 1 \otimes v^{(p)}), x_2) \\
= &\ {\rm Res}_{x_1} (x_1 - x_2)^{-1} Y(v^{(1)}\otimes \cdots
\otimes v^{(p-1)} \otimes 1, x_1)Y(1 \otimes \cdots
\otimes 1 \otimes v^{(p)}, x_2)\\
&\ - {\rm Res}_{x_1} (-x_2 + x_1)^{-1}Y(1 \otimes \cdots \otimes 1
\otimes v^{(p)}, x_2)Y(v^{(1)}\otimes \cdots \otimes v^{(p-1)}
\otimes 1, x_1). \end{split}\end{equation*} It follows that $(v^{(1)}
\otimes \cdots \otimes v^{(p)})_n$ is a
linear combination of the operators $(v^{(1)}\otimes \cdots \otimes
v^{(p-1)} \otimes 1)_{n_1} \cdot (1 \otimes \cdots \otimes 1 \otimes
v^{(p)})_{n_2}$. Thus the lemma holds by the
induction hypothesis.\epf\\

Now we define the notion of tensor product module for tensor product
strongly $A = A_1 \oplus \cdots \oplus A_p$-graded conformal vertex
algebra $V = V_1 \otimes \cdots \otimes V_p$ with the notions above.
Let $\tilde{A_1}, \dots, \tilde{A_p}$ be abelian groups containing
$A_1, \dots, A_p$ as subgroups, respectively, and let $W_1, \dots,
W_p$ be strongly $\tilde{A_1}, \dots, \tilde{A_p}$-graded modules
for $V_1, \dots, V_p$, respectively.\\

Let \[\tilde{A} = \tilde{A_1} \oplus \cdots \oplus \tilde{A_p},\]
Then we can construct the\ {\it tensor product strongly
$\tilde{A}$-graded module}
\[W = W_1\otimes \cdots \otimes W_p\] for the
tensor product strongly $A$-graded algebra $V$ by means of the
definition
\[ Y(v^{(1)}\otimes \cdots \otimes v^{(p)}, x) = Y(v^{(1)},x)\otimes \cdots
\otimes Y(v^{(p)}, x)\;\;\; \mbox{for}\; v^{(i)} \in V_i,\; i = 1, \dots, p,
\]
\[L(n) = L_1(n)\otimes 1 \otimes \cdots \otimes 1 + \cdots + 1\otimes
\cdots \otimes 1\otimes L_p(n)\;\;\; \mbox{for}\; n \in \mathbb{Z}.
\]
(Here we use the notation $L_i(n)$ for the operators associated with
$\omega^i$ on $W_i$, $i = 1, \dots,
p$.) The $\tilde{A}$-grading of $W$ is defined as \[W =
\coprod_{\beta \in \tilde{A}} W^{(\beta)},\] with \[ W^{(\beta)} =
 W_1^{(\beta_1)}\otimes
\cdots \otimes W_p^{(\beta_p)},
\]
where $\beta_i \in \tilde{A_i}$, $i = 1, \dots, p$, are such that
$\beta_1+ \cdots + \beta_p = \beta$. The $\mathbb{C}$-grading of $W$
is defined as \[W = \coprod_{n \in \mathbb{C}} W_{(n)},\] where
\[W_{(n)} = \sum_{n_1 + \cdots + n_p = n}(W_1)_{(n_1)}\otimes \cdots
\otimes (W_p)_{(n_p)}.\]
It follows that the $\mathbb{C}$-grading is given by the operator
$L(0)$ on $W$ defined above. It is clear that the algebra $V$ is also
a module for itself.\\

\begin{propo}\label{W}The structure $W$ constructed above is a
strongly $\tilde{A}$-graded module for the tensor product strongly
$A$-graded conformal vertex algebra $V$.
\end{propo}

\begin{assum}\label{countable group}{\rm In the remainder of this
paper, we always assume that $A$, and that each $A_i$ ($i = 1, \cdots,
p$) is a countable abelian group.}\end{assum}

Using Proposition \ref{x}, we now prove:

\begin{theo}\label{y}
Let $W = W_1 \otimes \cdots \otimes W_p$ be a strongly $\tilde{A} =
\tilde{A_1} \oplus \cdots \oplus \tilde{A_p}$-graded $V$-module,
with the notations as above. Then $W$ is irreducible if and only if
each $W_i$ is irreducible.
\end{theo}

$Proof$. The ``only if" part is trivial. For the ``if" part, for
simplicity of notation, we take $p = 2$ without losing any essential
content. Take a nonzero submodule $W \subset W_1\otimes W_2$, let
$w_1^{(1)}, \dots, w_n^{(1)} \in W_1$ and $w_1^{(2)}, \dots,
w_n^{(2)} \in W_2$ be linearly independent such that $\Sigma_{j =
1}^{n} a_j(w_j^{(1)} \otimes w_j^{(2)}) \in W$, where each $a_j \neq
0$. Take any $w^{(1)} \in W_1$, $w^{(2)} \in W_2$. By Proposition
\ref{x}, the commuting ring consists of
 the scalars for $W_1$ and $W_2$. Thus by the density theorem (see for example
Section 5.8 of [J]), there are $b_1 \in A(V_1; W_1 \otimes W_2),\
b_2 \in A(V_2; W_1 \otimes W_2)$ such that
\[b_1 \cdot w_1^{(1)} = w^{(1)},\ b_1\cdot w_i^{(1)} = 0, \ {\rm for}\ i = 2, \dots, n.\]
\[b_2\cdot w_1^{(2)} = w^{(2)},\ b_2\cdot w_i^{(2)} = 0, \ {\rm for}\ i = 2, \dots, n.\]
Then
\[ (b_1b_2)\cdot\Sigma_{j=1}^{n}\ a_j (w_j^{(1)}\otimes w_j^{(2)} ) = a_1 (w^{(1)}\otimes w^{(2)})\ \in W.\]
Hence $w^{(1)} \otimes w^{(2)} \in W$, and so $W = W_1 \otimes
W_2$. \epf \\

\section{Strongly $(\mathfrak{h}, A)$-graded vertex algebras and their
strongly $(\mathfrak{h}, \tilde{A})$-graded modules}

For some strongly $A$-graded vertex algebras $V$, there is a vector
space $\mathfrak{h}$ consisting of mutually commuting operators
induced by $V$ such that the $A$-grading of $V$ is given by
$\mathfrak{h}$ in the following way: for $\alpha \in A$,
$V^{(\alpha)}$ is the weight space of $\mathfrak{h}$ of weight
$\alpha$. Here is an example:

\begin{exam}\label{weight operator example}{\rm Consider the strongly $L$-graded conformal
vertex algebra $V_L$ in Example \ref{lattice vertex algebra}. For $h
\in \mathfrak{h}$, there is an operator $h(0)$ on $V_L$ such
that
\[h(0)\cdot V_L^{(\alpha)} = \langle h, \alpha\rangle
V_L^{(\alpha)}.\] We identify $\mathfrak{h}$ with the set of operators
\[\{ h(0)\ =\ (h(-1)\cdot {\bf 1})_0\ |\ h \in \mathfrak{h}\}\](see
Chapter 8 of \cite{FLM}).
Since the symmetric bilinear form $\langle \cdot, \cdot \rangle$ is
nondegenerate, $V_L^{(\alpha)}$ is characterized as the weight space
of $\mathfrak{h}$ of weight $\alpha$.
 }\end{exam}

Consider the tensor algebra $T(V[t, t^{-1}])$ over the vector space
$V[t, t^{-1}]$. Then any $V$-module $W$, in particular, $V$ itself,
can be regarded as a $T(V[t, t^{-1}])$-module uniquely determined by
the condition that for $v \in V$, $n \in \mathbb{Z}$, $v \otimes
t^n$ acts on $W$ as $v_n$. In the following definitions, we consider
a particular subspace of $T(V[t, t^{-1}])$ acting on $V$ and $W$.

\begin{defi}\label{strongly {h, A} graded algebra}{\rm A strongly
$A$-graded vertex algebra equipped with a vector subspace \[
\mathfrak{h} \subset T(V[t, t^{-1}]) \] is called {\it strongly
($\mathfrak{h}$, $A$)-graded} if there is a nondegenerate pairing
\begin{equation*}\begin{split}
\langle\cdot, \cdot\rangle: \mathfrak{h} \times A &\longrightarrow \mathbb{C}\\
(h, \alpha) &\longmapsto \langle h, \alpha\rangle
\end{split}\end{equation*}
linear in the first variable and additive in the second variable,
such that $\mathfrak{h}$ acts commutatively on $V$ and
\[V^{(\alpha)} = \{v \in V\ |\ h\cdot v = \langle h, \alpha\rangle v,
{\rm\ for\ all}\ h \in \mathfrak{h}\}.\]}\end{defi}

By Definition \ref{strongly {h, A} graded algebra}, the strongly
graded conformal vertex algebra $V_L$ in Example \ref{weight operator
example} is strongly $(\mathfrak{h}, L)$-graded, where $\mathfrak{h}$
is the set of operators $\{(h(-1)\cdot {\bf 1})_0\ |\ h \in
L\otimes_\mathbb{Z} \mathbb{C}\}$.

For a strongly ($\mathfrak{h}$, $A$)-graded vertex algebra $V$, a
natural module category is the category of strongly
$\tilde{A}$-graded $V$-modules $W$ with an action of $\mathfrak{h}$,
such that the $\tilde{A}$-grading on $W$ is given by weight spaces
of $\mathfrak{h}$. Here is an example:

\begin{exam}\label{strongly {h, A}-graded module example}{\rm As in Example \ref{lattice vertex
algebra}, any sublattice $M$ of $L^{\circ}$ containing $L$ gives rise
to a strongly $M$-graded $V_L$-module $V_M$. Take $\mathfrak{h} = L
\otimes_{\Z}\C$ and identify $\mathfrak{h}$ as the set of operators
$\{(h(-1)\cdot {\bf 1})_0\ |\ h \in \mathfrak{h}\}$ as in Example
\ref{weight operator example}. Then for $\beta \in M$, \[V_M^{(\beta)}
= \{w \in V_M\ |\ h\cdot w = \langle h, \beta\rangle w, {\rm\ for\
all}\ h \in \mathfrak{h}\}.\] so that we have examples of the
following:}\end{exam}

\begin{defi}\label{abelian group weight operator module}{\rm A
strongly $\tilde{A}$-graded module for a
strongly ($\mathfrak{h}$, $A$)-graded vertex algebra is said to be
{\it strongly ($\mathfrak{h}$, $\tilde{A}$)-graded} if
there is a nondegenerate pairing
\begin{equation*}\begin{split}
\langle\cdot, \cdot\rangle: \mathfrak{h} \times \tilde{A} &\longrightarrow \mathbb{C}\\
(h, \beta) &\longmapsto \langle h, \beta\rangle
\end{split}\end{equation*}
linear in the first variable and additive in the second variable,
such that the operators in $\mathfrak{h}$ act commutatively on $W$
and
\[W^{(\beta)} = \{w \in W\ |\ h\cdot w = \langle h, \beta\rangle w,
{\rm\ for\ all}\ h \in \mathfrak{h}\}.\]}\end{defi}

\begin{rema}\label{sub}{\rm Submodules and quotient modules of strongly $(\mathfrak{h},
\tilde{A})$-graded conformal modules are also strongly
$(\mathfrak{h}, \tilde{A})$-graded modules. Irreducible strongly
$(\mathfrak{h}, \tilde{A})$-graded modules are strongly
$(\mathfrak{h}, \tilde{A})$-graded modules without nontrivial
submodules. Strongly $(\mathfrak{h}, \tilde{A})$-graded module
homomorphisms are strongly $\tilde{A}$-graded module homomorphisms
which commute with the actions of $\mathfrak{h}$.}
\end{rema}

The following propositions are natural analogues of Proposition
\ref{V} and Proposition \ref{W}.

\begin{propo}\label{tensor product algebra}Let $V_1, \dots, V_p$ be strongly
$(\mathfrak{h}_1, A_1), \dots, (\mathfrak{h}_p, A_p)$-graded
conformal vertex algebras, respectively. Let $A = A_1 \oplus \cdots
\oplus A_p$, $\mathfrak{h} = \mathfrak{h}_1 \oplus \dots \oplus
\mathfrak{h}_p$, and let $\langle \cdot, \cdot\rangle_i$ denote the
pairing between $\mathfrak{h_i}$ and $A_i$, for $i = 1, \dots, p$.
Then the tensor product algebra $V = V_1 \otimes \cdots \otimes V_p$
becomes a strongly $(\mathfrak{h}, A)$-graded conformal vertex
algebra, where the nondegenerate pairing is given by:
\begin{equation*}\begin{split}
\langle\cdot, \cdot\rangle: \mathfrak{h} \times A &\longrightarrow \mathbb{C}\\
(h, \alpha) &\longmapsto \sum_{i = 1}^p \langle h_i,
\alpha_i\rangle_i,
\end{split}\end{equation*} where $h = h_1 + \dots + h_p$, $\alpha
= \alpha_1 + \dots + \alpha_p$,
for $h_i \in \mathfrak{h}_i$, $\alpha_i \in A_i$, $i = 1, \dots, p$,
and \[V^{(\alpha)} = V_1^{(\alpha_1)} \otimes \cdots \otimes
V_p^{(\alpha_p)} = \{v \in V_1 \otimes \cdots \otimes V_p \ | \ h
\cdot v = \langle h , \alpha \rangle v, {\rm\ for\ all}\ h \in
\mathfrak{h} \}.\]
\end{propo}

{\it Proof.} It is easy to see that the pairing defined above is
nondegenerate, and $V^{(\alpha)}$ is characterized uniquely as the
eigenspace of $\mathfrak{h}$. \epf

\begin{propo}\label{tensor product module}Let $W_1, \dots, W_p$ be
strongly $(\mathfrak{h}_1, \tilde{A_1}), \dots, (\mathfrak{h}_p,
\tilde{A_p})$-graded conformal modules for strongly
$(\mathfrak{h}_1, A_1), \dots, (\mathfrak{h}_p, A_p)$-graded
conformal vertex algebras $V_1, \dots, V_p$, respectively. Let
$\tilde{A} = \tilde{A_1} \oplus \cdots \oplus \tilde{A_p}$,
$\mathfrak{h} = \mathfrak{h}_1 \oplus \dots \oplus \mathfrak{h}_p$,
and let $\langle \cdot, \cdot\rangle_i$ denote the pairing between
$\mathfrak{h_i}$ and $\tilde{A_i}$, for $i = 1, \dots, p$. Then the
tensor product module $W = W_1 \otimes \cdots \otimes W_p$ becomes a
strongly $(\mathfrak{h}, \tilde{A})$-graded module for the
strongly graded vertex algebra $V$, where the nondegenerate pairing
is given by:
\begin{equation*}\begin{split}
\langle\cdot, \cdot\rangle: \mathfrak{h} \times \tilde{A} &\longrightarrow \mathbb{C}\\
(h, \beta) &\longmapsto \sum_{i = 1}^p \langle h_i,
\beta_i\rangle_i,
\end{split}\end{equation*} where $h = h_1 + \dots + h_p$, $\beta = \beta_1 + \dots + \beta_p$,
for $h_i \in \mathfrak{h}_i$, $\beta_i \in \tilde{A_i}$, $i = 1,
\dots, p$, and \[W^{(\beta)} = W_1^{(\beta_1)} \otimes \cdots
\otimes W_p^{(\beta_p)} = \{w \in W_1 \otimes \cdots \otimes W_p \ |
\ h \cdot w = \langle h , \beta \rangle w, {\rm\ for\ all}\ h \in
\mathfrak{h} \}.\ \ \ \  \square\]
\end{propo}

The following proposition is an analogue and consequence of Theorem \ref{y}.

\begin{theo}\label{y1}Let $W = W_1 \otimes \cdots \otimes W_p$ be a
strongly $(\mathfrak{h}, \tilde{A})$-graded module constructed in
Proposition \ref{tensor product module}. Then $W$ is irreducible if
and only if each $W_i$ is irreducible. \end{theo}

\section{Irreducible modules for tensor product strongly graded algebra}
Our goal is to determine all the strongly $(\mathfrak{h},
\tilde{A})$-graded irreducible modules for the tensor product
strongly $(\mathfrak{h}, A)$-graded conformal vertex algebra
constructed in Proposition \ref{tensor product algebra}. To do this,
we need to define a more specific kind of strongly $(\mathfrak{h},
\tilde{A})$-graded modules as follows:

\begin{defi}\label{separate}{\rm Let $V_1, \dots, V_p, V$ be strongly $(\mathfrak{h}_1, A_1), \dots,
(\mathfrak{h}_p, A_p), (\mathfrak{h}, A)$-graded conformal
vertex algebras, respectively, as in the setting of Proposition
\ref{tensor product algebra}. Let $W$ be a strongly $(\mathfrak{h},
\tilde{A})$-graded $V$-module, where $\tilde{A}$ is an abelian
group containing $A$ as a subgroup, so that in particular, for
$\beta \in \tilde{A}$,
\[W^{(\beta)} = \{w \in W\ |\ h\cdot w = \langle h, \beta\rangle w,
{\rm\ for\ all}\ h \in \mathfrak{h}\}.\]
Assume that there exists an abelian subgroup $\tilde{A_i}$ of
$\tilde{A}$ containing $A_i$ as a subgroup for each $i = 1, \dots,
p$ such that
\begin{equation*}\begin{split}
&\tilde{A} = \tilde{A_1} \oplus \cdots \oplus \tilde{A_p},                 \\
&\langle \mathfrak{h}_i, \tilde{A_j}\rangle = 0\   {\rm if}\ i \neq
j
\end{split}\end{equation*}
and such that $W$ is a doubly graded $V_i$-module with respect to
$\tilde{A_i}$ and the $\tilde{A_i}$-grading is given by
$\mathfrak{h}_i$ in the following way: For $\beta_i \in
\tilde{A_i}$,
\[W^{(\beta_i)} = \{w \in W\ |\ h_i\cdot w = \langle h_i, \beta_i
\rangle w, {\rm\ for\ all}\ h_i \in \mathfrak{h}_i\}.\] Then $W$ is
called a {\it strongly} $((\mathfrak{h}_1, \tilde{A_1}), \dots,
(\mathfrak{h}_p, \tilde{A_p}))$-{\it graded
$V$-module}.}\end{defi}

\begin{rema}\label{sub1}{\rm Submodules and quotient modules
of strongly $((\mathfrak{h}_1, \tilde{A_1}), \dots,
(\mathfrak{h}_p, \tilde{A_p}))$-graded
$V$-modules are also strongly
$((\mathfrak{h}_1, \tilde{A_1}), \dots,
(\mathfrak{h}_p, \tilde{A_p}))$-graded modules. Irreducible strongly
$((\mathfrak{h}_1, \tilde{A_1}), \dots,
(\mathfrak{h}_p, \tilde{A_p}))$-graded modules are strongly
$((\mathfrak{h}_1, \tilde{A_1}), \dots,
(\mathfrak{h}_p, \tilde{A_p}))$-graded modules without nontrivial
submodules. Strongly $((\mathfrak{h}_1, \tilde{A_1}), \dots,
(\mathfrak{h}_p, \tilde{A_p}))$-graded module
homomorphisms are strongly $(\mathfrak{h}, \tilde{A})$-graded $V$-module homomorphisms.}
\end{rema}

\begin{exam}\label{strongly example}{\rm The strongly $(\mathfrak{h},
\tilde{A})$-graded tensor product module $W_1 \otimes \cdots \otimes
W_p$ constructed in Proposition
\ref{tensor product module} is a strongly $((\mathfrak{h}_1,
\tilde{A_1}), \dots, (\mathfrak{h}_p, \tilde{A_p}))$-graded $V_1
\otimes \cdots \otimes V_p$-module.}\end{exam}

{}From Example \ref{lattice vertex algebra}, we can see that any
$V_L$-module is a strongly $L^{\circ}$-graded module. Based on this
fact, it is easy to show that the following example satisfies the
conditions in Definition \ref{separate}.

\begin{exam}\label{strongly example1}{\rm Let $V^{\natural}$ be the
moonshine module constructed in \cite{FLM}, which is a strongly
$(\langle 0 \rangle, \langle 0 \rangle)$-graded conformal vertex
algebra as in Example \ref{vertex operator algebra}; let $V_L$ be the
conformal vertex algebra associated with the even 2-dimensional
unimodular Lorentzian lattice $L$, which is a strongly
$(\mathfrak{h}, L)$-graded conformal vertex algebra as constructed
in Example \ref{lattice vertex algebra}. Then any strongly
$(\mathfrak{h}, L)$-graded module for $V^{\natural}\otimes V_L$ is
strongly $((\langle 0 \rangle, \langle 0 \rangle), (\mathfrak{h},
L))$-graded (note that $L$ is a self-dual lattice, i.e., $L^{\circ} =
L$).}\end{exam}

\begin{nota}\label{multi2}{\rm For $\beta_1 \in \tilde{A_1}, \dots,
\beta_p \in \tilde{A_p}$, we let $W^{(\beta_1, \dots, \beta_p)}$
denote the following common weight space of $\mathfrak{h}_1, \dots,
\mathfrak{h}_p$, i.e.,
\[W^{(\beta_1, \dots, \beta_p)}: = \{w\in W\ |\ h_i\cdot w = \langle
h_i, \beta_i\rangle w,\ {\rm for\ all}\ h_i \in \mathfrak{h}_i, \ i
= 1, \dots, p \}.\]}
\end{nota}

Next we assume $W$ to be a strongly $((\mathfrak{h}_1,
\tilde{A_1})$, \dots, $(\mathfrak{h}_p, \tilde{A_p}))$-graded
$V_1 \otimes \cdots \otimes V_p$-module, with the notation as in
Definition \ref{separate}.

\begin{propo}\label{multi3}Suppose that $W$ is irreducible. Then for $\beta_1 \in \tilde{A_1},
\dots, \beta_p \in \tilde{A_p}$, $W^{(\beta_1, \dots, \beta_p)}$ is
irreducible under the algebra of operators $A(V_1 \otimes \cdots
\otimes V_p; W^{(\beta_1, \dots, \beta_p)})$.
\end{propo}

{\em Proof.} The proof is similar to the proof of Proposition
\ref{c}. \epf

\begin{lemma}\label{separate1}For $\beta \in \tilde{A}$,
we have\[W^{(\beta)} = W^{(\beta_1, \dots,
\beta_p)},\] where $\beta = \beta_1 + \cdots + \beta_p$.\end{lemma}

{\em Proof.} This is a consequence of Definition
\ref{separate}. \epf

\begin{theo}\label{z}Let $W$ be a strongly
$((\mathfrak{h}_1, \tilde{A_1})$, \dots, $(\mathfrak{h}_p,
\tilde{A_p}))$-graded irreducible $V_1\otimes \cdots \otimes
V_p$-module, with the notions as in Definition \ref{separate}. Then
$W$ is a tensor product of irreducible
strongly $(\mathfrak{h}_i, \tilde{A_i})$-graded $V_i$-modules, for $i
=1, \dots, p$.
\end{theo}

$Proof$. For simplicity of notation, we take $p = 2$, as above.
Since $W$ is irreducible, by Remark \ref{congruent}, $W =
\coprod_{\bar{n}= \mu} W_{(n)}$ for some $\mu \in
\mathbb{C}/\mathbb{Z}$, where $\bar{n}$ denotes the equivalent class
of $n\in \mathbb{C}$ in $\mathbb{C}/\mathbb{Z}$. Choose $\beta \in \tilde{A}$
such that $W^{(\beta)} \neq 0$. Then there exists $n_0 \in \mathbb{C}$ such that
$W_{(n_0)}^{(\beta)}$ is the lowest weight space of $W^{(\beta)}$.
Since $W_{(n_0)}^{(\beta)}$ is finite dimensional and we are working
over $\mathbb{C}$, there exists a simultaneous eigenvector $w_0 \in
W_{(n_0)}^{(\beta)}$ for the commuting operators $L_i(0)$ and the
operators in $\mathfrak{h}_i,\ i = 1, 2$. Denote by $n_1, n_2 \in
\mathbb{Z}$ the corresponding eigenvalues for $L_1(0),\ L_2(0)$.
Then we have $n_0 = n_1 + n_2$. Denote by $\beta_1 \in \tilde{A_1},\
\beta_2 \in \tilde{A_2}$ the corresponding weights for
$\mathfrak{h}_1,\ \mathfrak{h}_2$. By Lemma \ref{separate1},
we have $W^{(\beta)} = W^{(\beta_1, \beta_2)}$,
and $\beta = \beta_1 + \beta_2$.\\

Now the $L(-1)$-derivative condition and the $L(0)$-bracket formula
imply that \[ [L_1(0), Y(v^{(1)}\otimes 1, x)] = Y(L_1(0)(v^{(1)}
\otimes 1), x) + x\frac{d}{dx}Y(v^{(1)} \otimes 1, x)
\]
for $v^{(1)} \in V_1$. Thus for doubly homogeneous vector $v^{(1)}$
and $n \in \mathbb{Z}$,
\[
\mathrm{wt}_1(v^{(1)} \otimes 1)_n = \mathrm{wt}_1 (v^{(1)}\otimes
1)- n- 1,
\]
where wt$_1$ refers to $L_1(0)$-eigenvalue on both $V_1\otimes V_2$
and the space of operators on $W$. In particular, $(v^{(1)} \otimes
1)_n$ permutes $L_1(0)$-eigenspaces. Moreover, since $(1 \otimes
v^{(2)})_n$, for $v^{(2)} \in V_2$, commutes with $L_1(0)$, it preserves
$L_1(0)$-eigenspaces. Of course, similar statements hold for
$L_2(0)$, $\mathfrak{h}_1(0)$, $\mathfrak{h}_2(0)$.\\

By Lemma \ref{multi3}, $W^{(\beta_1, \beta_2)}$ is irreducible under
the algebra of the operators $A(V_1 \otimes V_2; W^{(\beta_1,
\beta_2)})$. Then $W^{(\beta_1, \beta_2)}$ is generated by $w_0$ by
the irreducibility, and is spanned by elements of the form \[
(v_1^{(1)}\otimes 1)_{m_1} \cdots (v_k^{(1)}\otimes 1)_{m_k}(1
\otimes v_1^{(2)})_{n_1}\cdots (1 \otimes v_l^{(2)})_{n_l}w_0
\] where $v_i^{(1)} \in V_1$ and $v_j^{(2)} \in V_2$, $v_i^{(1)},
v_j^{(2)}$ are doubly homogeneous, and the $A$-weights of $\sum_{i=1}^m\
v_i^{(1)}$ and $\sum_{j=1}^n
\ v_j^{(2)}$ are 0.\\

Hence $W^{(\beta_1, \beta_2)}$ is the direct sum of its simultaneous
eigenspaces for $L_i(0)$ and $\mathfrak{h}_i$, for $i = 1, 2$, and
the $L_1(0)$, $L_2(0)$-eigenvalues are bounded below by $n_1$,
$n_2$, respectively. It follows that the lowest weight space
$W_{(n_0)}^{(\beta_1, \beta_2)}$ is filled up by the simultaneous
eigenspace for the operators $L_i(0)$ with eigenvalues $n_i$. To be
more precise, we use $W_{(n_1, n_2)}^{(\beta_1, \beta_2)}$ to denote
the subspace $W_{(n_0)}^{(\beta_1,
\beta_2)}$. By a similar argument as in Proposition \ref{multi3}, $W_{(n_1,
n_2)}^{(\beta_1, \beta_2)}$ is irreducible under the algebra of
operators $A(V_1 \otimes V_2; W_{(n_1,
n_2)}^{(\beta_1, \beta_2)})$.\\

By the density theorem, the algebra $A(V_1 \otimes V_2; W_{(n_1,
n_2)}^{(\beta_1, \beta_2)})$ fills up End $W_{(n_1, n_2)}^{(\beta_1,
\beta_2)}$. Because  $A(V_1; W_{(n_1, n_2)}^{(\beta_1, \beta_2)})$
and $A(V_2; W_{(n_1, n_2)}^{(\beta_1, \beta_2)})$ are commuting
algebras of operators and $A(V_1 \otimes V_2; W_{(n_1,
n_2)}^{(\beta_1, \beta_2)})$ is generated by  $A(V_1; W_{(n_1,
n_2)}^{(\beta_1, \beta_2)})$ and $A(V_2; W_{(n_1, n_2)}^{(\beta_1,
\beta_2)})$, we see that \[ {\rm End}\ W_{(n_1, n_2)}^{(\beta_1,
\beta_2)} = A(V_1; W_{(n_1, n_2)}^{(\beta_1, \beta_2)})A(V_2;
W_{(n_1, n_2)}^{(\beta_1, \beta_2)}).\] Choose an irreducible
$A(V_1; W_{(n_1, n_2)}^{(\beta_1, \beta_2)})$-submodule $M_1$ of
$W_{(n_1, n_2)}^{(\beta_1, \beta_2)}$. Then $A(V_1; W_{(n_1,
n_2)}^{(\beta_1, \beta_2)})$ acts faithfully on $M_1$ since any
element of $A(V_1; W_{(n_1, n_2)}^{(\beta_1, \beta_2)})$ which
annihilates $M_1$ annihilates\\ $A(V_2; W_{(n_1, n_2)}^{(\beta_1,
\beta_2)})\cdot M_1\ =\ A(V_2; W_{(n_1, n_2)}^{(\beta_1,
\beta_2)})A(V_1; W_{(n_1, n_2)}^{(\beta_1, \beta_2)})\cdot M_1\ =\
(\mathrm{End}\ W_{(n_1, n_2)}^{(\beta_1, \beta_2)})M_1\ =\ W_{(n_1,
n_2)}^{(\beta_1, \beta_2)}$. Thus $A(V_1; W_{(n_1, n_2)}^{(\beta_1,
\beta_2)})$ restricts faithfully to End $M_1$ and hence is
isomorphic to a full matrix algebra. Similarly, $A(V_2; W_{(n_1,
n_2)}^{(\beta_1, \beta_2)})$ is isomorphic to a full matrix algebra.
It follows that
\[ \mathrm{End}\ W_{(n_1, n_2)}^{(\beta_1, \beta_2)} = A(V_1;
W_{(n_1, n_2)}^{(\beta_1, \beta_2)}) \otimes A(V_2; W_{(n_1,
n_2)}^{(\beta_1, \beta_2)}).
\] Then $W_{(n_1, n_2)}^{(\beta_1, \beta_2)}$ has the structure
\[ W_{(n_1, n_2)}^{(\beta_1, \beta_2)} = M_1 \otimes M_2\]
as an irreducible $A(V_1; W_{(n_1, n_2)}^{(\beta_1, \beta_2)})
\otimes A(V_2; W_{(n_1, n_2)}^{(\beta_1, \beta_2)})$-module. Here,
as an irreducible $A(V_i; W_{(n_1, n_2)}^{(\beta_1,
\beta_2)})$-submodule of $W^{(\beta_1, \beta_2)}_{(n_1, n_2)}$,
$M_i$ has $\tilde{A_i}$-grading $\beta_i$ induced by
$\mathfrak{h}_i$, and has $\mathbb{C}$-grading $n_i$ induced by
$L_i(0)$, respectively, for $i\ =\ 1,2$.\\

Let \[ w^0 = y_1 \otimes y_2\] (where $y_i \in M_i$, for $i = 1, 2$)
be a nonzero decomposable tensor in $W_{(n_1, n_2)}^{(\beta_1,
\beta_2)}$. Let $W_i$ be the doubly graded $V_i$-submodule of $W$
generated by $w^0$. Then the module $W_1$ has a strongly
$(\mathfrak{h}_1,
\tilde{A_1})$-graded $V_1$-module structure such that
\[W_1 = \coprod_{n \in \mathbb{C},\ \gamma \in \tilde{A}}
(W_1)_{(n)}^{(\gamma)},
\] where \begin{eqnarray*}(W_1)_{(n)}^{(\gamma)} = {\rm
span}\{(v_1^{(1)}\otimes {\bf 1})_{s_1}\cdots (v_p^{(1)}\otimes{\bf
1})_{s_p}w^0|{\rm wt}\ v_1^{(1)} - s_1 - 1 + \cdots + {\rm wt}\
v_p^{(1)} - s_p -1 = n - n_1,\\ A\mbox{-}{\rm wt}\ v_1^{(1)} +
\cdots + A\mbox{-}{\rm wt}\ v_p^{(1)} = \gamma - \beta_1,\
v_1^{(1)}, \dots, v_p^{(1)}\in V_1,\ s_1, \dots, s_p \in
\mathbb{Z}\}.\end{eqnarray*}

This module we constructed satisfies the grading restrictions
(\ref{set:dmltc}) and (\ref{set:dmfin}) in Definition
\ref{def:dgw}, which follows from the fact that $W$ is a strongly
graded $V_1 \otimes V_2$-module and each doubly homogeneous subspace
of $W_1$ lies in the doubly homogeneous subspace of $W$. Also,
$W_1^{(\gamma)}$ is the weight space of $\mathfrak{h}_1$ with weight
$\gamma$, hence by Definition \ref{abelian group weight operator
module},
$W_1$ is a strongly $(\mathfrak{h}_1, \tilde{A_1})$-graded $V_1$-module.\\

We claim that $W_1$ is $V_1$-irreducible (and similarly for $W_2$).
In fact, consideration of the abelian group grading shows that any
nonzero $V_1$-submodule of $W_1$ not intersecting $W^{(\beta_1,
\beta_2)}$ will give rise to a nonzero $V_1 \otimes V_2$-submodule
of $W$ not intersecting $W^{(\beta_1, \beta_2)}$. Thus any nonzero
$V_1$-submodule of $W_1$ must intersect $W^{(\beta_1, \beta_2)}$.
 Then consideration of the weight shows that the $(\beta_1, \beta_2)$-subspace of any nonzero
$V_1$-submodule of $W_1$ not intersecting $W_{(n_1, n_2)}^{(\beta_1,
\beta_2)}$ would give rise to a nonzero $A(V_1 \otimes V_2;
W^{(\beta_1, \beta_2)})$-submodule of $W^{(\beta_1, \beta_2)}$ not
intersecting $W_{(n_1, n_2)}^{(\beta_1, \beta_2)}$. Thus any nonzero
$V_1$-submodule of $W_1$ must intersect $W_{(n_1, n_2)}^{(\beta_1,
\beta_2)}$. But the irreducible $A(V_1; W_{(n_1, n_2)}^{(\beta_1,
\beta_2)})$-module $A(V_1; W_{(n_1, n_2)}^{(\beta_1, \beta_2)})\cdot
w^0$ is the full intersection of $W_1$ and $W_{(n_1,
n_2)}^{(\beta_1, \beta_2)}$, so that the $V_1$-submodule must
contain $w^0$ and hence be all of
$W_1$. This proves the $V_1$-irreducibility of $W_1$.\\

Finally, to show that $W$ is isomorphic to $W_1\otimes W_2$,
consider the abstract tensor product $V_1\otimes V_2$-module
$W_1\otimes W_2$, where $W_i$ is the strongly $\tilde{A_i}$-graded
$V_i$-module defined above, for $i = 1, 2$. Define a linear map
\[
\varphi: W_1\otimes W_2 \rightarrow W
\]
\[
b_1\cdot w^0 \otimes b_2 \cdot w^0 \mapsto b_1b_2\cdot w^0,
\]
where $b_i$ is any operator induced by $V_i$. Then $\varphi$ is well
defined and is a $V_1 \otimes V_2$-module homomorphism. Since $W_1
\otimes W_2$ is irreducible by Theorem \ref{y}, $\varphi$ is a
module isomorphism.\epf

\begin{exam}\label{lattice case}{\rm Let $V_{L_i}$ be the conformal vertex algebra associated
with an even lattice $L_i$ as in Example \ref{lattice vertex algebra},
where $i = 1, \dots, p$. Let $V_{L_1} \otimes \cdots \otimes V_{L_p}$
be the tensor product strongly graded vertex algebra of $V_{L_1},
\dots, V_{L_p}$. By the construction of a lattice vertex algebra in
Example \ref{lattice vertex algebra}, we have \[V_{L_1} \otimes \cdots
\otimes V_{L_p} = V_{L_1 \oplus \cdots \oplus L_p},\]and every
irreducible $V_{L_1 \oplus \cdots \oplus L_p}$-module is equivalent to
a module of the form \[V_{L_1 + \gamma_1 \oplus \cdots \oplus L_p +
\gamma_p} = V_{L_1 + \gamma_1} \otimes \cdots \otimes V_{L_p +
\gamma_p},\] for some $\gamma_i \in L_i^{\circ}$, $i = 1, \dots, p$.
This example illustrates Theorem \ref{z}.}\end{exam}

Now we can describe our main examples:

\begin{corol}\label{main theorem for simplicity}The only irreducible strongly
$(\mathfrak{h}, L)$-graded module of $V^\natural \otimes V_L$,
where $L$ is the unique even 2-dimensional unimodular Lorentzian
lattice and $\mathfrak{h} = \{(h(-1) \cdot 1)_{0}\ |\ h \in L
\otimes_{\mathbb{Z}} \mathbb{C}\}$, up to equivalence, is
itself.\end{corol}

$Proof$. Let $W$ be an irreducible strongly $(\mathfrak{h},
L)$-graded module of $V^\natural \otimes V_L$. Then by Example
\ref{strongly example1}, $W$ is a strongly $((\langle 0 \rangle,
\langle 0 \rangle), (\mathfrak{h}, L))$-graded module of
$V^\natural \otimes V_L$. By Theorem \ref{z}, it is a tensor product
of an irreducible strongly $(\langle 0 \rangle, \langle 0
\rangle)$-graded $V^\natural$-module with an irreducible strongly
$(\mathfrak{h}, L)$-graded $V_L$-module. By [D2], $V^\natural$ is
its only irreducible module, up to equivalence. Also, by [D1]
(cf. \cite{LL}, Example \ref{lattice vertex algebra}), $V_L$
is its only irreducible module because $L$ is self-dual. Therefore
\[ W\ =\ V^\natural \otimes V_L
\]as claimed. \epf

\begin{rema}{\rm In Corollary \ref{main theorem for simplicity}, the
2-dimensional self-dual Lorentzian lattice can of course be
generalized to any self-dual nondegenerate even lattice.}
\end{rema}

\section{Complete reducibility}
\begin{defi}\label{ssm}{\rm Let $V$ be a strongly $(\mathfrak{h},
A)$-graded conformal vertex algebra. Then a strongly
$(\mathfrak{h}, \tilde{A})$-graded $V$-module is called {\it
completely reducible} if it is a direct sum of irreducible strongly
$(\mathfrak{h}, \tilde{A})$-graded $V$-modules.}
\end{defi}

\begin{nota}{\rm In the remainder of this section, we will always let $A =
A_1 \oplus \cdots \oplus A_p$, $\mathfrak{h} = \mathfrak{h}_1 \oplus \cdots \oplus \mathfrak{h}_p$,
and $V = V_1 \otimes \cdots \otimes V_p$.}\end{nota}

\begin{defi}\label{SS}{\rm A strongly $((\mathfrak{h}_1, \tilde{A_1}), \dots, (\mathfrak{h}_p,
\tilde{A_p}))$-graded module for the tensor product conformal
vertex algebra $V$ is called {\it completely reducible} if it is a
direct sum of irreducible strongly $((\mathfrak{h}_1,
\tilde{A_1}), \dots, (\mathfrak{h}_p, \tilde{A_p}))$-graded
$V$-modules.}
\end{defi}

\begin{theo}\label{TSS} Let $V_1, \dots, V_p$ be strongly $(\mathfrak{h}_1, A_1), \dots,
(\mathfrak{h}_p, A_p)$-graded conformal vertex algebras,
respectively, and let $V$ be their tensor product strongly $(
\mathfrak{h}, A)$-graded conformal vertex algebra. Then every
strongly $((\mathfrak{h}_1, \tilde{A_1}), \dots, (\mathfrak{h}_p,
\tilde{A_p}))$-graded $V$-module is completely reducible if and
only if every strongly $(\mathfrak{h}_i, \tilde{A_i})$-graded
$V_i$-module is completely reducible.
\end{theo}
{\it Proof.} It suffices to prove the result for $n = 2$. Let $W$ be
a strongly $((\mathfrak{h}_1, \tilde{A_1}), (\mathfrak{h}_2,
\tilde{A_2}))$-graded $V = V_1 \otimes V_2$-module. Then by
Proposition \ref{separate1}, we can take $w \in W^{(\beta_1,
\beta_2)}_{(n_1, n_2)}$, where $\beta_i \in \tilde{A_i}$, $n_i \in
\mathbb{C}$, for $i = 1, 2$. \\

Let $M$ be the strongly $((\mathfrak{h}_1, \tilde{A_1}), (\mathfrak{h}_2,
\tilde{A_2}))$-graded $V_1\otimes V_2$-submodule of $W$ generated
by $w$, i.e., $M$ is spanned by elements of the form
\begin{eqnarray*}
(v_1^{(1)}\otimes {\bf 1})_{s_1}\cdots
(v_p^{(1)}\otimes{\bf 1})_{s_p}({\bf 1} \otimes v_1^{(2)})_{t_1}
\cdots ({\bf 1} \otimes v_q^{(2)})_{t_q}w
\end{eqnarray*}
where
$v_1^{(1)}, \dots,
v_p^{(1)}$ are doubly homogeneous elements in $V_1$ and $v_1^{(2)}, \dots,
v_q^{(2)}$ are doubly homogeneous elements in $V_2$, respectively, and $s_1, \dots,
s_p, t_1, \dots, t_q \in \mathbb{Z}$. Let $M_i$ be the doubly graded
$V_i$-submodule of $M$ generated by $w$.
Then $M_i$ is a strongly $(\mathfrak{h}_i, \tilde{A_i})$-graded $V_i$-module,
respectively, for $i = 1, 2$, in an obvious way as in the proof of
Theorem \ref{z}.\\

By Proposition \ref{tensor product module} and Example \ref{strongly example},
$M_1\otimes M_2$ is strongly $((\mathfrak{h}_1, \tilde{A_1}),
(\mathfrak{h}_2, \tilde{A_2}))$-graded. Moreover, we have a
strongly $((\mathfrak{h}_1, \tilde{A_1}), (\mathfrak{h}_2,
\tilde{A_2}))$-graded $V_1\otimes V_2$-module epimorphism from
 $M_1\otimes M_2$ to $M$ by sending $b_1w \otimes b_2w \mapsto b_1b_2w$, where
$b_i$ is an operator induced by $V_i$, for $i = 1,\ 2$. If every
strongly $(\mathfrak{h}_i, \tilde{A_i})$-graded $V_i$-module is
completely reducible, then $M_i$ is a direct sum of irreducible
strongly $(\mathfrak{h}_i, \tilde{A_i})$-graded $V_i$-modules and
therefore $M_1\otimes M_2$ is a direct sum of irreducible strongly
$((\mathfrak{h}_1, \tilde{A_1}), (\mathfrak{h}_2,
\tilde{A_2}))$-graded $V_1\otimes V_2$-modules (see Theorem
\ref{y1}). Then as a quotient module of $M_1\otimes M_2$, $M$ is
also a direct sum of irreducible strongly $((\mathfrak{h}_1,
\tilde{A_1}), (\mathfrak{h}_2, \tilde{A_2}))$-graded $V_1\otimes
V_2$-modules, and consequently, $W$ is a direct sum of irreducible
strongly $((\mathfrak{h}_1, \tilde{A_1}), (\mathfrak{h}_2,
\tilde{A_2}))$-graded
$V_1\otimes V_2$-modules.\\

Conversely, assume that every strongly $((\mathfrak{h}_1,
\tilde{A_1}), (\mathfrak{h}_2, \tilde{A_2}))$-graded $V_1\otimes
V_2$-module $W$ is completely reducible. We first observe that $V_1
\otimes V_2$ is strongly $((\mathfrak{h}_1, A_1),
(\mathfrak{h}_2, A_2))$-graded, hence a $((\mathfrak{h}_1,
\tilde{A_1}), (\mathfrak{h}_2, \tilde{A_2}))$-graded $V_1\otimes
V_2$-module itself by Proposition \ref{tensor product algebra} and
Example \ref{strongly example}, and hence is a direct sum of irreducible
strongly $((\mathfrak{h}_1, \tilde{A_1}), (\mathfrak{h}_2,
\tilde{A_2}))$-graded modules. Let $W$ be an irreducible strongly
$((\mathfrak{h}_1, \tilde{A_1}), (\mathfrak{h}_2,
\tilde{A_2}))$-graded $V_1 \otimes V_2$-module. Then $W$ is a
tensor product of an irreducible strongly $(\mathfrak{h}_1,
\tilde{A_1})$-graded module for $V_1$ and an irreducible strongly
$(\mathfrak{h}_2, \tilde{A_2})$-graded module for $V_2$ by Theorem
\ref{z}. In particular, $V_1$ has irreducible strongly
$(\mathfrak{h}_1, \tilde{A_1})$-graded modules and  $V_2$ has irreducible strongly
$(\mathfrak{h}_2, \tilde{A_2})$-graded modules,
respectively.\\

Let $W_1$ be a strongly $(\mathfrak{h}_1, \tilde{A_1})$-graded
$V_1$-module and $W_2$ be an irreducible strongly $(\mathfrak{h}_2,
\tilde{A_2})$-graded $V_2$-module. Since every strongly
$((\mathfrak{h}_1, \tilde{A_1}), (\mathfrak{h}_2,
\tilde{A_2}))$-graded $V_1\otimes V_2$-module is completely
reducible, $W_1 \otimes W_2$ is a direct sum of irreducible strongly
$((\mathfrak{h}_1, \tilde{A_1}), (\mathfrak{h}_2,
\tilde{A_2}))$-graded modules:\[W_1\otimes W_2 = \coprod_i M_i\]
where each $M_i$ is an irreducible strongly $((\mathfrak{h}_1,
\tilde{A_1}), (\mathfrak{h}_2, \tilde{A_2}))$-graded $V_1\otimes
V_2$-module. Fix $i$ and let $x_1^{(i)}, \dots, x_n^{(i)} \in W_1$
and $y_1^{(i)}, \dots, y_n^{(i)} \in W_2$ be linearly independent
doubly homogeneous elements such that $\sum_j c_jx_j^{(i)} \otimes
y_j^{(i)} \in M_i$, where $c_j \in \mathbb{C}, c_j \neq 0$. By the
density theorem (as in the proof of Theorem \ref{y}), each
$x_j^{(i)} \otimes y_j^{(i)} \in M_i$. Let $W_{i1}$ be the doubly graded $V_1$-submodule
of $W_1$ generated by $x_{j_0}^{(i)}$, for some $j_0 \in \{1, 2 , \dots, n\}$. Then $W_{i1}$
is a strongly $(\mathfrak{h}_1, \tilde{A_1})$-graded $V_1$-submodule as in the
proof of Theorem \ref{z}. By the irreducibility of $M_i$, we see that
$M_i = W_{i1} \otimes W_2$ and that $W_{i1}$ is an irreducible
strongly $(\mathfrak{h}_1, \tilde{A_1})$-graded $V_1$-submodule of
$W_1$. Therefore, $W_1\otimes W_2 = (\coprod_i W_{i1})\otimes W_2$.
By the density theorem, for any nonzero $w_2 \in W_2$, $W_1 \otimes w_2
= (\coprod_i W_{i1})\otimes w_2$. Hence as a $V_1$-module, $W_1 \cong
(\coprod_i W_{i1})$, and thus $W_1$ is completely reducible. Similarly for
$V_2$. \epf

\begin{exam}{\rm Let $V_{L_i}$ be the conformal vertex algebra associated
with an even lattice $L_i$ as in Example \ref{lattice vertex algebra},
where $i = 1, \dots, p$. Let $V_{L_1} \otimes \cdots \otimes V_{L_p}$
be the tensor product strongly graded vertex algebra of $V_{L_1},
\dots, V_{L_p}$. By the construction of a lattice vertex algebra as in
Example \ref{lattice vertex algebra}, we have \[V_{L_1} \otimes \cdots
\otimes V_{L_p} = V_{L_1 \oplus \cdots \oplus L_p}.\]As in Example
\ref{lattice vertex algebra}, every module for $V_{L_1 \oplus \cdots
\oplus L_p}$, hence for $V_{L_1} \otimes \cdots \otimes V_{L_p}$, is
completely reducible. This example illustrates Theorem \ref{TSS}.}
\end{exam}

\begin{corol}Every strongly $(\mathfrak{h}, L)$-graded module for the strongly $(\mathfrak{h}, L)$-graded conformal vertex algebra $V^\natural \otimes V_L$,
where $L$ is the unique even 2-dimensional unimodular Lorentzian
lattice and $\mathfrak{h} = \{(h(-1) \cdot 1)_{0}\ |\ h \in L
\otimes_{\mathbb{Z}} \mathbb{C}\}$, is completely reducible.
\end{corol}

\def\refname{\hfil{REFERENCES}}

\noindent {\small \sc Department of Mathematics, Rutgers University,
110 Frelinghuysen Rd., Piscataway, NJ 08854-8019}

\vspace{1em}

\noindent {\em E-mail address}: yookinwi@math.rutgers.edu

\end{document}